\documentclass[12pt]{article}

\usepackage{amsmath, amssymb}

\newtheorem{thm}{Theorem}
\newtheorem{prop}{Proposition}

\title{Toric  3-folds defined by quadratic binomials}
\author{Shoetsu Ogata}

\begin{document}
\maketitle

\begin{abstract}
Let $(X,A)$ be a polarized nonsingular toric 3-fold 
with $\Gamma(X, K_X+A)=0$.
Then for any ample line bundle $L$ on $X$ the image of the
embedding $\Phi_L: X\to \mathbb{P}(\Gamma(X,L))$ is an intersection of quadrics.
\end{abstract}
\section*{Introduction}

Sturmfels asked in \cite{S2} whether a nonsingular projective
toric variety should be
defined by only quadrics if it is embedded by global sections of
a normally generated ample line bundle.
An evidence  has been obtained by Koelman
(\cite{K3}) before Sturmfels asked the question.
Koelman showed that projective toric surfaces are defined by
binomials (differences of two monomials) of degree at most three
(\cite{K2}) and obtained a criterion when the
surface needs  defining equations of degree three (\cite{K3}).
He used combinatorics of plane polygons.

Let $X$ be a projective algebraic variety and let $L$ an ample line
bundle on it.  If the natural homomorphism
\begin{equation}\label{0;0}
\phi: 
S:=\mbox{Sym}\ \Gamma(X, L) \longrightarrow R:=\bigoplus_{k\ge0}\Gamma(X, L^{\otimes k})
\end{equation}
is surjective, then
Mumford (\cite{Mf}) calls $L$ {\it normally generated}.
A normally generated ample line bundle is always  very ample, but not conversely.
We define the ideal $I(X,L)$ of the graded ring $S$ as
$$
I(X,L):=\mbox{Ker}\ \phi \subset S=\bigoplus_{k\ge0}S_k.
$$
Mumford (\cite{Mf}) calls $L$ {\it normally presented} if $I(X,L)$ is generated by elements of degree two.

Let $X$ be a toric variety of dimension $n$ and $L$  an ample line bundle on it.
In general, $L$ is not very ample if $n\ge3$.  On the other hand, $L^{\otimes k}$ is normally
generated for $k\ge n-1$(\cite{EW}), and the ideal $I(X, L^{\otimes k})$ is generated by quadrics for $k\ge n$ (\cite{BGT},\cite{NO}), or for $k\ge n-1$ and $n\ge3$ (\cite{Og0}).

We  know that there exists a polarized toric variety $(X, L)$ of dimension
$n\ge3$ such that $L$ is very ample but $L^{\otimes (n-2)}$ is not normally generated(\cite{BG}, \cite{Og2}).
We also know that any ample line bundle on a {\it nonsingular} toric variety is
always very ample (\cite{D}, see also \cite[Corollary 2.15]{Od}).
Ogata (\cite{Og1}) showed that an ample line bundle $L$ on a nonsingular toric
3-fold $X$  is normally generated if 
 the adjoint bundle $L+K_X$ is not big.

In this paper we give a partial answer to Sturmfels' question.

\begin{thm}\label{tm1}
Assume that $(X,A)$ is a polarized nonsingular toric variety of dimension three with
$\Gamma(X, A+K_X)=0$.  Let $L$ be an ample line bundle on $X$ and 
$\Phi_L: X\to \mathbb{P}(\Gamma(X,L))$ the associated embedding.
Then the image $\Phi_L(X)$ is the common zero of quadratic binomials.
\end{thm}

The proof is separated into two propositions  as Proposition \ref{p3:1} in Section 3 and Proposition \ref{p5:1} in Section 5.

In Section 1 we recall the basic fact about toric varieties and ample line bundles on them, and corresponding
lattice polytopes.
In Section 2 we give an algerbo-geometric proof of the result of Koelman and explain the classification of
$(X,A)$ satisfying the condition in Theorem~\ref{tm1}.
In Section 3 we discuss the binomials defining affine parts of $\Phi_L(X)$ and give a strategy to prove
Theorem~\ref{tm1}.
In Section 4 we point out some property of nonsingular lattice polygons  (Proposition \ref{p4:1}).
In Section 5 we give a proof of the main part of Theorem~\ref{tm1} as Proposition \ref{p5:1}.

\bigskip

\section{Polarized toric varieties}

In this section we recall the fact about toric varieties and ample line bundles on them
and corresponding lattice polytopes (see, for example, Oda's book \cite{Od} or Fulton's book \cite{Fu}).

Let $M$ be a free abelian group of rank $n$ and $M_{\Bbb R}:=
M\otimes_{\Bbb Z}\Bbb R \cong \Bbb R^n$ the extension of coefficients.
Set $\mathbb{C}[M]$ the group algebra of $M$ and $T:=\mbox{Spec}\mathbb{C}[M]
\cong (\mathbb{C}^{\times})^n$ the algebraic torus of dimension $n$.
Then the group of characters $\text{Hom}_{\text{gr}}(T, \Bbb C^{\times})$ is isomorphic to $M$.
For an element $m\in M$ we denote by $\bold e(m): T\to \mathbb{C}^{\times}$  the character corresponding to $m$. 

A toric variety $X$ is a normal algebraic variety with an algebraic
action $T\times X \to X$
of the algebraic torus $T$ such that $X$ contains an open orbit $O$ isomorphic to $T$ and that
the action is compatible with the inclusion $T\cong O\to X$ and the multiplication $T\times T\to T$.

We define a {\it lattice polytope} as the convex hull $P:= \mbox{Conv}\{m_1,
\dots, m_r\}$ of a finite subset $\{m_1, \dots, m_r\}$ of $M$ in $M_{\mathbb{R}}$.
We define the dimension of a lattice polytope $P$ as that of the smallest
affine subspace $\mathbb{R}(P)$ containing $P$.

Let $X$ be a projective toric variety of dimension $n$ and $L$ an ample line bundle
on $X$.   Then there exists a lattice polytope $P$ of dimension $n$ such that
the space of global sections of $L$ is described by
\begin{equation}\label{0:1}
\Gamma(X, L)\cong \bigoplus_{m\in P\cap M}\mathbb{C}{\bf e}(m),
\end{equation}
where ${\bf e}(m)$ is considered as a rational function on $X$ since
$T$ is identified with the dense open subset (see
 \cite[Section 2.2]{Od} or \cite[Section 3.5]{Fu}).
 We also have
 \begin{equation}\label{0:2}
 \Gamma(X, L\otimes \omega_X)\cong \bigoplus_{m\in \mbox{\small{int}}(P)\cap M}\mathbb{C}{\bf e}(m),
\end{equation}
where $\omega_X$ is the dualizing sheaf of $X$.
 
Conversely, for a lattice polytope $P$ in $M_{\mathbb{R}}$ of dimension $n$
set $V(P)$ the set of all vertices of $P$.
For each vertex $v\in V(P)$ define the convex cone $C_v(P):=\mathbb{R}_{\ge0}(P-v)$ and
the affine toric variety $U_v:=\mbox{Spec}\ \mathbb{C}[C_v(P)\cap M]$.
We obtain an toric variety by gluing them:
$$
X=\bigcup_{v\in V(P)} U_v.
$$
We define a line bundle $L$ so that
$$
\Gamma(U_v,L)=\bold e(v) \mathbb{C}[C_v(P)\cap M].
$$
Then $L$ is ample and  satisfies the equality (\ref{0:1})
(see \cite[Chaper 2]{Od} or \cite[Section 1.5]{Fu}).

Let $A$ and $B$ be two ample line bundles on $X$, and $P_A$ and $P_B$ the corresponding
lattice polytopes.  Then $A\otimes B$ corresponds to the Minkowski sum
$$
P_A+P_B:=\{x+y\in M_{\mathbb{R}}: \ x\in P_A \ \mbox{and}\ y\in P_B\}
$$
(see \cite[Section 1.5]{Fu}).

If $X$ is nonsingular, then all $U_v$ are isomorphic to $\mathbb{C}^n$.  This implies that
there exists a $\mathbb{Z}$-basis $\{m_1, \dots, m_n\}$ of $M$ such that
$$
C_v(P)=\mathbb{R}_{\ge0}m_1+\dots +\mathbb{R}_{\ge0}m_n
$$
(see \cite[Theorem 1.10]{Od}).

\bigskip

\section{Algebro-geometric approach}

We recall the results of Koelman.  He treated the case of dimension two.

\begin{thm}[\cite{K2}, \cite{K3}]\label{tk1}
Any ample line bundle $L$ on a projective toric surface $X$ is normally generated and 
the ideal $I(X,L)$ is generated by elements of degree at most three.
Moreover, it is generated by quadrics unless $\Gamma(X, L\otimes \omega_X)\not=0$ and
$\dim\Gamma(X,L)-\dim\Gamma(X, L\otimes \omega_X)=3$.
\end{thm}

In his proof Koelman uses combinatrics of lattice polygons.
Let $P$ be the lattice polygon corresponding to a polarized toric surface $(X,L)$.
The conditions in the exception are $\mbox{int}(P)\cap M\not=\emptyset$
and the number of lattice points in the boundary $\partial P$ of $P$ is equal to three.
Thus $P$ is a triangle and $X$ is a singular toric surface isomorphic to $\mathbb{P}^2/G$.

Here we give a proof of Theorem \ref{tk1} by using a method of projective algebraic geometry.
Let $C\in |L|$ be a general member of the linear system of $L$.  Then
$C$ is a nonsingular curve of genus $g=\sharp(\mbox{int}(P)\cap M)$.
Let $L_C$ denote the restriction to $C$.  Then we have
$$
\deg L_C =\sharp(\partial P\cap M) +2g-2.
$$
Since $P$ is a convex polygon, $\sharp(\partial P\cap M)\ge3$.  The theorem of Fujita (\cite{Fj})
says that $L_C$ is normally generated if $deg L_C\ge 2g+1$ and that $I(c, L_C)$ is generated
by quadrics if $\deg L_c\ge 2g+2$.
By regular ladder theorem (\cite{Fj}), we see that $L$ is always normally generated, and that 
$I(X,L)$ is generated by only quadrics if $\partial P\cap M)\ge4$.

Next, we consider the case of dimension three.
Ogata (\cite{Og1}) classified the polarized toric 3-folds satisfying the condition in Theorem \ref{tm1}.
\begin{prop}[\cite{Og1}]\label{p2:1}
Let $(X,A)$ be a nonsingular polarized toric variety of dimension three with $\Gamma(X,A+K_X)=0$.
Then $X$ is one of the followings.
\begin{enumerate}
\item[\rm (1)] a blow up $\mathbb{P}^3$ along at most 4 invariant points,
\item[\rm (2)] a blow up $\mathbb{P}^2$-bundle over $\mathbb{P}^1$ along at most 2 invariant points,
\item[\rm(3)] a $\mathbb{P}^1$-bundle over a nonsingular toric surface.
\end{enumerate}
\end{prop}

Let $M=\mathbb{Z}^3$ with a basis $\{e_1, e_2, e_3\}$.
Let $Q$ be the lattice polytope of dimension three corresponding $(X,A)$ in Proposition \ref{p2:1}.
The condition $\Gamma(X, A+K_X)=0$ implies that $\mbox{int}(Q)\cap M=\emptyset$.
Set $\Delta_3:=\mbox{Conv}\{0, e_1, e_2, e_3\}$ the basic 3-simplex.

In the case (1), $Q$ is $k\Delta_3$ for $1\le k\le3$, one cut of $\Delta_3$ from $2\Delta_3$ or
at most 4 cuts of $\Delta_3$'s from $3\Delta_3$.  See Figure 1.

In the case (2), $Q$ is a prism with the base $\Delta_2=\mbox{Conv}\{0, e_1,e_2\}$ and three edges of
length $a, b, c\ge1$, or at most one cut from the base and the roof of a prism with the base $2\Delta_2$ and three edges of length $d, e, f\ge1$
such that $e-f$ and $e-d$ are both even.  See Figure 2.

In the case (3), $Q$ has parallel two facets $F_0$ and $F_1$ width one such that $F_i$ is a lattice polygon 
corresponding to a polarized nonsingular toric surface $(Y, L_i)$.  $F_0$ and $F_1$ have the same number of edges and corresponding edges are parallel.

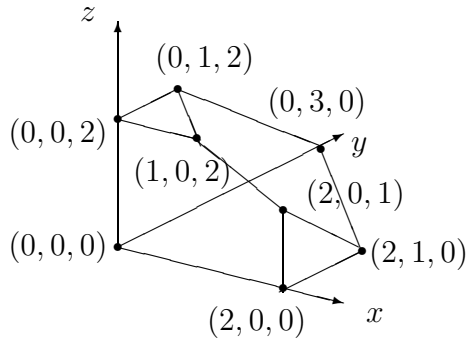
\begin{figure}[h]
 \begin{center}
 \setlength{\unitlength}{1mm}
  \begin{picture}(50,40)(10,10)
   \put(20,20){\vector(4,-1){30}}
   \put(20,20){\vector(0,1){30}}
   \put(20,20){\vector(2,1){30}}
   \put(2,15){\makebox(20,10){$(0,0,0)$}}
   \put(53,10){\makebox(10,10)[bl]{$x$}}
   \put(15,50){\makebox(10,10)[bl]{$z$}}
   \put(51,32){\makebox(10,10)[bl]{$y$}}
   \put(41.5,14){\line(2,1){11}}
   \put(42,14){\line(0,1){11}}
   \put(28,41){\line(5,-2){19}}
   \put(28,41){\line(2,-5){3}}
   \put(20,37){\line(4,-1){11}}
   \put(20,37){\line(2,1){8}}
   \put(47,33){\line(2,-5){5.5}}
   \put(42,25){\line(2,-1){10}}
   \put(42,25){\line(-5,4){12}}
   \put(35,5){\makebox(10,10)[r]{$(2,0,0)$}}
   \put(55,14){\makebox(10,10)[]{$(2,1,0)$}}
   \put(45,22){\makebox(10,10)[l]{$(2,0,1)$}}
   \put(43,37){\makebox(10,10)[br]{$(0,3,0)$}}
   \put(42,25){\circle*{1}}
   \put(47,33){\circle*{1}}
   \put(25,40){\makebox(10,10)[l]{$(0,1,2)$}}
   \put(22,22){\makebox(10,10)[tl]{$(1,0,2)$}}
   \put(2,30){\makebox(20,10){$(0,0,2)$}}
   \put(20,20){\circle*{1}}
   \put(42,14.5){\circle*{1}}
   \put(20,37){\circle*{1}}
   \put(52.5,19.5){\circle*{1}}
   \put(30.5,34.5){\circle*{1}}
   \put(28,41){\circle*{1}}
  \end{picture}
  \end{center}
  \caption{typical $Q$ of (1)}
\end{figure}

 \begin{figure}[h]
 \begin{center}
 \begin{tabular}{lr}
  \setlength{\unitlength}{1mm}
  \begin{picture}(50,40)(15,15)
   \put(20,20){\vector(4,-1){30}}
   \put(20,20){\vector(0,1){30}}
   \put(20,20){\vector(2,1){30}}
   \put(2,15){\makebox(20,10){$(0,0,0)$}}
   \put(53,10){\makebox(10,10)[bl]{$x$}}
   \put(15,50){\makebox(10,10)[bl]{$z$}}
   \put(51,32){\makebox(10,10)[bl]{$y$}}
   \put(40,15){\line(0,1){30}}
   \put(35,27){\line(0,1){23}}
   \put(35,28){\line(1,-3){4.5}}
   \put(20,35){\line(1,1){15}}
   \put(20,35){\line(2,1){20}}
   \put(35,50){\line(1,-1){5}}
   \put(35,6){\makebox(10,10){$(1,0,0)$}}
   \put(40,20){\makebox(10,10){$(0,1,0)$}}
   \put(43,40){\makebox(10,10){$(1,0,a)$}}
   \put(30,47){\makebox(10,10){$(0,1,b)$}}
   \put(5,30){\makebox(10,10){$(0,0,c)$}}
   \put(20,20){\circle*{1}}
   \put(40,15){\circle*{1}}
   \put(35,27){\circle*{1}}
   \put(35,50){\circle*{1}}
   \put(20,35){\circle*{1}}
   \put(40,45){\circle*{1}}
  \end{picture}&

  \setlength{\unitlength}{1mm}
  \begin{picture}(50,40)(5,15)
   \put(20,20){\vector(4,-1){30}}
   \put(20,20){\vector(0,1){30}}
   \put(20,20){\vector(2,1){30}}
   \put(7,15){\makebox(10,10){$(0,0,0)$}}
   \put(53,10){\makebox(10,10)[bl]{$x$}}
   \put(15,50){\makebox(10,10)[bl]{$z$}}
   \put(51,32){\makebox(10,10)[bl]{$y$}}
   \put(40,24){\line(0,1){21}}
   \put(35,34){\line(0,1){16}}
   \put(20,35){\line(1,1){15}}
   \put(20,35){\line(2,1){20}}
   \put(35,50){\line(1,-1){5}}
   \put(37,30){\makebox(10,10){$(0,2,1)$}}
   \put(42,18){\makebox(10,10){$(2,0,1)$}}
   \put(46,40){\makebox(10,10){$(2,0,d+1)$}}
   \put(30,48){\makebox(10,10){$(0,2,e+1)$}}
   \put(3,30){\makebox(10,10){$(0,0,f-1)$}}
   \put(20,20){\circle*{1}}
   \put(35,50){\circle*{1}}
   \put(20,35){\circle*{1}}
   \put(40,45){\circle*{1}}
   \put(20,20){\circle*{1}}
   \put(40,24){\circle*{1}}
  \put(35,34){\circle*{1}}
   \put(20,35){\circle*{1}}
   \put(40,45){\circle*{1}}
   \put(30,17.5){\circle*{1}}
   \put(28,24){\circle*{1}}
   \put(35,34){\line(1,-2){5}}
    \put(28,24){\line(1,-3){2.3}}
    \put(30,17.5){\line(3,2){10}}
    \put(28,24){\line(2,3){7}}
     \put(25,9){\makebox(10,10){$(1,0,0)$}}
     \put(16,21){\makebox(10,10){$(0,1,0)$}}
   
  \end{picture}
  \end{tabular}
 \end{center}
\caption{typical $Q$ of (2)}
 \label{fig0}
\end{figure}
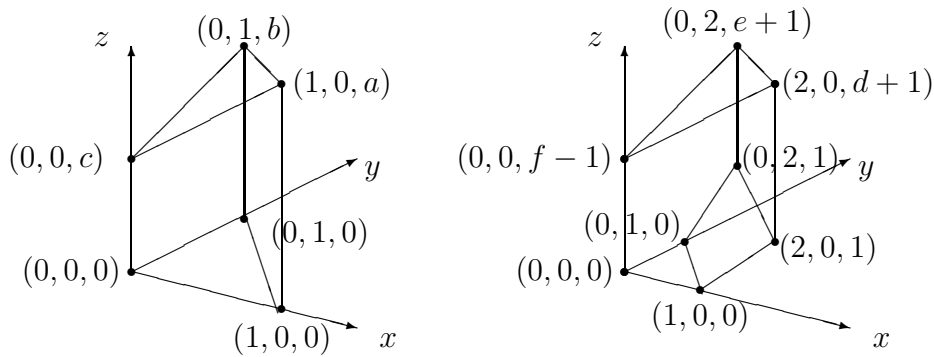

\begin{prop}\label{p2:2}
Let $(X,A)$ be a polarized nonsingular toric 3-fold inProposition \ref{p2:1}.
If $X$ is (1) or (2), then the ideal $I(X,A)$ is generated by only quadrics.
\end{prop}

{\it Proof}.  Let $S_1, S_2\in |A|$ be two general members of the linear system of $A$.
Set $C=S_1\cap S_2$. Then $C$ is a nonsingular curve. Set $g(C)$ the genus of $C$.
We have
\begin{eqnarray*}
\chi(A_C) &=&\sharp(Q\cap M) -\sharp(\mbox{int}(Q)\cap M)=\sharp(Q\cap M),\\
g(C) &=& \sharp(\mbox{int}(2Q)\cap M) -2\sharp(\mbox{int}(Q)\cap M)=\sharp(\mbox{int}(2Q)\cap M).
\end{eqnarray*}
From Riemann-Roch formula we have
$$
\deg A_C=g(C) -1 +\chi(A_C).
$$
If $\chi(A_C)\ge g(C)+3$, Fujita's theorem (\cite{Fj}) says that $I(X,A)$ is generated by quadrics.
By calculation of the numbers $\sharp(Q\cap M)$ and $\sharp(\mbox{int}(2Q)\cap M)$,
we obtain a proof.   \hfill  $\Box$

\bigskip

\section{Ideals of affine parts}

When a projective variety $X$ is embedded by the complete linear system of a very ample line bundle $L$,
$\Phi_L: X\to \mathbb{P}^r$, it is an effective way to investigate polynomials defining 
its affine parts $\Phi_L(X)\cap (\mathbb{P}^r\setminus H_i)$ in order to find polynomials defining $\Phi_L(X)$.
We know that toric varieties are defined by binomials (\cite{ES}).

Let $P\subset M_{\mathbb{R}}$ be a lattice polytope corresponding to a nonsingular polarized toric variety
$(X,L)$ of dimension three.
Since $X$ is a union of affine toric varieties $U_v$ ($v\in V(P)$), the embedding
$\Phi_L: X \to \mathbb{P}(\Gamma(X,L))$ is defined by that of  affine parts $U_v\subset X$:
$$
\Phi_L|U_v: U_v=\mbox{Spec}\ \mathbb{C}[C_v(P)\cap M] \to \mbox{Spec}(\mbox{Sym}\langle(P-v)\cap M\rangle_{\mathbb{C}})
$$
for all $v\in V(P)$.

Set $P\cap M=\{m_0, m_1, \dots, m_r\}$.  We assume that $m_0$ is a vertex of $P$ and that
$m_1,m_2,m_3$ are the  lattice points nearest to $m_0$ on three edges meeting each other at $m_0$.
Then the convex cone $C_{m_0}(P)$ is
$$
C_{m_0}(P)= \mathbb{R}_{\ge0}(m_1-m_0) +\mathbb{R}_{\ge0}(m_2-m_0) +\mathbb{R}_{\ge0}(m_3-m_0).
$$
Since $\{m_1-m_0, m_2-m_0, m_3-m_0\}$ is a $\mathbb{Z}$-basis of $M\cong \mathbb{Z}^3$,
$$
U_{m_0}=\mbox{Spec}\ \mathbb{C}[C_{m_0}(P)\cap M] \cong \mathbb{C}^3.
$$

Let $Z_0, Z_1, \dots, Z_r$ be the homogeneous coodinats of $\mathbb{P}(\Gamma(X,L))\cong \mathbb{P}^r$
corresponding to $P\cap M$.  We consider the affine part $\Phi_L(X)\cap(Z_0\not=0)$.
Set $x_i=Z_i/Z_0$.  Then $(Z_0\not=0)=\mbox{Spec}\ \mathbb{C}[x_1, \dots, x_r]\cong \mathbb{C}^r$.
Since $\{m_1-m_0, m_2-m_0, m_3-m_0\}$ is a $\mathbb{Z}$-basis of $M\cong \mathbb{Z}^3$,
for $i\ge4$ we can uniquely write as
$$
m_i-m_0=\sum_{j=1}^3 a_{ij}(m_j-m_0) \quad (a_{ij}\ge0).
$$
From this exprresion we define binomials as 
$$
f_i=x_i-\prod_{j=1}^3 x_j^{a_{ij}}.
$$
Then we have
$$
\mathbb{C}[C_{m_0}(P)\cap M] \cong \mathbb{C}[x_1, \dots, x_r]/(f_4, \dots, f_r).
$$

Here we define a property ``2-$D(m_0)$" :
For $m_i\in P\cap M\setminus \{m_0,m_1, m_2, m_3\}$ there exist $m_k, m_l\in
P\cap M\setminus\{m_0\}$ such that $m_i+m_0=m_k+m_l$, that is, $m_i-m_0=(m_k-m_0)+(m_l-m_0)$.

If $P$ satisfies the property 2-$D(m_0)$, then we define new binomials as
$$
g_i=x_i-x_kx_l,
$$
and we have equality of ideals $(f_4, \dots, f_r) =(g_4, \dots, g_r)$.
From $g_i$ we obtain homogeneous binomials $G_i:=Z_0Z_i-Z_kZ_l$, and we see that
the affine part $\Phi_L(X)\cap(Z_0\not=0)$ is the common zero set of $G_4, \dots, G_r$.

\begin{prop}\label{p3:1}
Let $P\subset M_{\mathbb{R}}$ be a lattice polytope corresponding to a nonsingular polarized toric variety
$(X,L)$ of dimension three.  Assume that $X$ is one of (1) and (2) in Proposition \ref{p2:1}.
Then for each vertex $v\in V(P)$, $P$ satisfies the property 2-$D(v)$. 
\end{prop}
{\it Proof}.
When $X$ is (1),  $P$ is at most 4 cuts of $l_i\Delta_3$ from $k\Delta_3$ such that $k\ge1$, $l_i\ge0$ and
$l_i+l_j<k$ for $i\not=j$.
If $m_0\in V(P)$ is a vertex of $k\Delta$, then we take a coordinates of $M$ as
$$
m_0=0, m_1=(1,0,0), m_2=(0,1,0), m_3=(0,0,1).
$$
Set $m_i=(a,b,c)$.  Then $a,b,c\ge0$.  When $c=0$, we have $a\ge1$ or $b\ge1$, hence, $(a-1,b,0)\in P$
or $(a,b-1,0)\in P$ and $(a,b,0)=(a-1,b,0)+(1,0,0)$ or $(a,b,0)=(a,b-1,0)+(0,1,0)$.
When $c\ge1$, we have $(a,b,c-1)\in P$ and $(a,b,c)=(a,b,c-1)+(0,0,1)$.

If $m_0$ is a vertex arising after cut of $l\Delta_3$ with $l\ge1$,  then we take a coordinates of $M$ as
$$
m_0=0, m_1=(1,0,-1), m_2=(0,1,-1), m_3=(0,0,1).
$$
Set $m_i=(a,b,c)$.  Then $a,b\ge0$.  When $c\ge0$, we can do the same procedure as above.
When $c<0$, we have $(a-1,b,c+1)\in P$ or $(a,b-1,c+1)\in P$.

When $X$ is (2), $P$ is at most one cut of $l_0\Delta_3$ from the base and that of $l_1\Delta_3$ from the roof of a prism with the base $k\Delta_2$ and three edges of length $d, e, f\ge1$
such that $e-f$ and $e-d$ are in $k\mathbb{Z}$ and $k\ge1$, $k>l_i\ge0$. 
If we take a coordinates as in (1), then we have a proof.   \hfill  $\Box$

\bigskip

\section{Nonsingular lattice polygons}

In order to prove Theorem\ref{tm1}, we have to treat the case that $X$ is a $\mathbb{P}^1$-bundle over a
nonsingular toric surface $Y$.

The lattice polytope $Q$ corresponding to (3) in Proposition \ref{p2:1} has two parallel facets $F_0$ and $F_1$
width one.  In order to compare lattice points on $F_0$ and $F_1$ we need to know  some information near opposite vertices.
Let $M'=\mathbb{Z}^2$. We call a lattice parallelogram $S\subset M'_{\mathbb{R}}$ to be a {\it basic diamond}
if $\sharp(S\cap M')=4$.

Let $F\subset M'_{\mathbb{R}}$ be a nonsingular lattice polygon with $s+1$ edges.
Let $u_0, u_1, \dots, u_s$ be vertices of $F$ numbered as counter-clockwise.
By an affine transform of $M'$, we may set as 
$$
u_0=0, \quad u_1=(a,0), \quad u_s=(0,b)
$$
with $a,b\ge1$.
Set $E_0=[u_0, u_1], E_s=[u_0, u_s]$ two edges of $F$ meeting at $u_0$.
If $\mbox{int}(F)\cap M'\not=\emptyset$, then the point $(1,1)$ is contained in the interior of $F$.

\begin{prop}\label{p4:1}
Let $F\subset M'_{\mathbb{R}}$ be a nonsingular lattice polygon with $s+1$ vertices
$u_0, u_1, \dots, u_s$ as above.
Assume $\mbox{int}(F)\cap M'\not=\emptyset$.
\begin{enumerate}
\item[\rm(1)] If $F$ has an edge $[u_{t-1}, u_t]$  parallel to $E_s$,
there exists a basic diamond $S$ contained in $F$ such that $u_t$ is a vertex of $S$, $[u_{t-1}, u_t]$ contains one edge of $S$ and that $S$ stays in $F$ after the vertex $m'$ of $S$ opposite to $u_t$ is transformed to the origin,
that is, $S-m'\subset F$.
\item[\rm(2)] When $F$ has no edges parallel to $E_0$ nor $E_s$, set $u_t$ the farthest vertex of $F$
from $u_0$.  Let $S\subset F$ be a basic diamond such that $u_t$ is a vertex of $S$ and $S$ has two edges contained in $[u_{t-1}, u_t]$ and $[u_t, u_{t+1}]$, respectively.  Set $m'\in S$ the vertex opposite to $u_t$.
Then $S-m'\subset F$. 
\end{enumerate}
\end{prop}
{\it Proof.}
First, consider the case (2).   Set $u_t=(p,q)$.  Then $p,q\ge1$.
Set $u_{t-1}=(p-k\alpha, q-k\beta), u_{t+1}=(p-l\gamma, q-l\delta)$ with $k,l\ge1$.  Since $u_t$ is the farthest from $u_0$ and $F$ has no edges parallel to $E_0$ nor $E_s$, we have $\alpha, \beta, \gamma, \delta\ge1$.  Since
$F$ is nonsingular, $\beta\gamma-\alpha\delta=1$.  Since $\mbox{int}(F)\cap M'\not=\emptyset$,
$(p-\alpha-\gamma, q-\beta-\delta)\in \mbox{int}(F)$.
Set
$$
S:=\mbox{Conv}\{(p-\alpha,q-\beta),(p,q),(p-\gamma, q-\delta), (p-\alpha-\gamma, q-\beta-\delta)\}
$$
and $m'=(p-\alpha-\gamma, q-\beta-\delta)\in \mbox{int}(F)$.  Then $S$ is a basic diamond and $S\subset F$.
From the convexity of $F$ we see $S-m'\subset F$.

When the case (1), since the edge $[u_{t-1},u_t]$ is parallel to $E_s$, we see $\alpha=0, \beta=\gamma=1$.
If $\delta\ge0$, then set 
$$
S:=\mbox{Conv}\{(p,q-1),(p,q),(p-1, q-\delta), (p-1, q-\delta-1)\}
$$
and $m'=(p-1, q-\delta-1)\in \mbox{int}(F)$.

If $\delta<0$, then set $S:=\mbox{Conv}\{(p,q-1),(p,q),(p-1,q),(p-1,q-1)\}$ and $m'=(p-1,q-1)$.
Then $m'\in \mbox{int}(F)$.
In both cases, $S-m'\subset F$ from convexity of $F$.
\hfill  $\Box$

\bigskip

\section{Proof of Theorem \ref{tm1}}

From the argument in Section 3, it is enough to prove the following proposition in order to 
obtain a proof of Theorem \ref{tm1}.

\begin{prop}\label{p5:1}
Let $P\subset M_{\mathbb{R}}$ be a lattice polytope corresponding to a nonsingular polarized toric variety
$(X,L)$ of dimension three.  Assume that $X$ is (3) in Proposition \ref{p2:1}.
Then for each vertex $v\in V(P)$, $P$ satisfies the property 2-$D(v)$. 
\end{prop}
{\it Proof}.
First, we consider the lattice polytope $Q$ corresponding to $(X,A)$.
$Q$ has two parallel facets $F_0$ and $F_1$ width one.
We may assume that a vertex $m_0$ of $Q$ is a vertex of $F_0$.
From an affine transform of $M$, we may set as $m_0$ is the origin and choose a basis $\{e_1, e_2, e_3\}$ of $M$
so that $e_1$ and $e_2$ are contained in edges of $F_0$ and $e_3$ is a vertex of $F_1$.
Set $M'=\mathbb{Z}e_1+\mathbb{Z}e_2$.  Then $M=M'\oplus \mathbb{Z}e_3$.
We may consider as $F_0, F_1\subset M'_{\mathbb{R}}$ and $Q=\mbox{Conv}\{ F_0\times 0, F_1\times e_1\}$.
Both of $F_0$ and $F_1$ have $s+1$ edges with $s\ge2$ and contain $e_1$ and $e_2$ in their edges.

Take $m_i\in Q\cap M\setminus\{m_0, e_1,e_2,e_3\}$.  Set $m_i=(a,b,c)$.  Then $a,b\ge0$ and 
$a\ge1$, $b\ge1$, $c=0$, or $c=1$.

When $s=2$, $(a-1,b,c)\in Q$ or $(a,b-1,c)\in Q$.  Hence $(a,b,c)=(a-1,b,c)+e_1$ or $(a,b,c)=(a,b-1,c);e_2$.

When $s=3$, both of $F_0$ and $F_1$ have at least one pair of parallel edges.  Assume that they have
 edges parallel to $[0,e_2]$.  Then $(a,b-1,c)\in Q$ if $b\ge1$.
 
 Set $s\ge4$.  Since $F_0$ and $F_1$ are nonsingular, they contain  lattice points in their interiors.
 We apply Proposition \ref{p4:1} to $F_0$.  We have the basic diamond $S\subset F_0$ and the lattice point
 $m'\in S$.  Set $\Bar{S}=S-m'$.  Then $S=\Bar{S}+m'$.  Since each edge of $F_1$ is parallel to corresponding edge of $F_0$,  the  basic diamond of $F_1$ is a parallel transform of $\Bar{S}$, that is, $\Bar{S}+m''$.
 Set $\Bar{S}=\{0, u_1', u_2', u_1'+u_2'\}$.
 
Consider the case  $m_i\in F_0\times0$.  Set $R_0:=\mbox{Conv}\{\Bar{S}, \Bar{S}+m'\}$.
 If $m_i\in R_0\setminus\{u_1', u_2'\}$, then there exists an $m_j\in R_0\cap M'\setminus\{0\}$ satisfying
 $m_i=m_j+u_1'$ or $m_i=m_j+u_2'$.
 When $m_i\notin R_0$, if it is contained in the side of $e_2$, then $m_i-e_2\in F_0$, if it is 
 contained in the side of $e_1$, then $m_i-e_1\in F_0$.
 After several steps, it moves in $R_0$.
 
 When $m_i\in F_1\times e_3$, set $R_1:=\mbox{Conv}\{\Bar{S}, \Bar{S}+m''\}$ and $m_i=\Bar{m}_i\times e_3$.
 If $\bar{m}_i=e_1$, then $m_i=e_1+e_3$.
We may set $\bar{m}_i\not= e_1, e_2$.
If $\bar{m}_i\in R_1\setminus\{u_1', u_2'\}$, then there exists an $\bar{m}_j\in R_1\cap M'\setminus\{0\}$ satisfying
 $\bar{m}_i=\bar{m}_j+u_1'$ or $\bar{m}_i=\bar{m}_j+u_2'$.
 Then we have
 $$
 m_i=\bar{m}_i\times e_3= \bar{m}_j\times e_3+u_1'\times0,\quad\mbox{or} \quad
 m_i= \bar{m}_j\times e_3+u_2'\times0.
 $$
The same method holds even if $\bar{m}_i\notin R_1$.

Next, for general $P$, we know that $P$ has also two parallel facets $F_0$ and $F_l$ width $l\ge1$
and each slice $F_k$ parallel to $F_0$ with width $k$ ($1\le k\le l$) is also a nonsingular lattice polygon
with the same number of edges parallel to corresponding edges of $F_0$.
 Since $m_i\in P\cap M\setminus\{m_0, m_1,m_2,m_3\}$ is contained one $F_k$,
 we can employ the same process by replacing $F_1$ with $F_k$.\hfill  $\Box$

\end{document}